# On A Certain Krull Symmetry Of a Noetherian Ring

C.L. Wangneo

#A-1 Jammu University.

New Campus Jammu.J&k -180006, India

## Abstract

In this note we show that if a noetherian ring R is left and right Krull-homogenous in which.

$$X = \{ P_i \in \mathrm{spec}.R \, / \, |R/P_i|_r = |R|_r \}$$

$$V = \{ Q_j \in \mathrm{spec}.R \, | \, |R/Q_j|_l = |R|_l \}$$

and $P = \cap_{P_i \in X} P_i$ and $Q = \cap_{Q_j \in V} Q_j$

then the following hold true;

(1) $P^n = 0$, for some integer $n \geq 1$

(2) $X = V$.

Introduction:- The main theorem of this paper, namely theorem(4) concerns a certain Krull symmetry of a noetherian ring which is both left and right Krull homogenous. In particular we show that if R is a noetherian left and right Krull homogenous ring in which:

$$X = \{ P_i \; \varepsilon \; \text{spec}.R \; / \; |R/P_i|_r = |R|_r \}$$

$$V = \{ Q_j \; \varepsilon \; \text{spec}.R \; | \; |R/Q_j|_l = |R|_l \}$$

and $P = \cap_{P_i \varepsilon X} P_i$ and $Q = \cap_{Q_j \varepsilon V} Q_j$

then the following hold true;

(1) $P^n = 0$, for some integer $n \geq 1$

(2) $X = V$

Terminology:- For a few words about the terminology in this paper we mention that all our rings are with identity element and all our modules are unitary. A ring R is noetherian means that R is a left as well as right noetherian ring. For a right R module M the definitions of right Krull dimension, the associated prime ideals of M and the critical composition series of M we refer the reader to [4]. M is said to be $\alpha$ Krull-homogenous if every non zero sub module of M has Krull dimension $\alpha$. For the definition of a ring satisfying the large condition we refer the reader to[1] and for the definition of weak ideal invariance(w.i.i for short) of an ideal we refer the reader to [2]. We now present the list of symbols that we shall use throughout:

If R is a ring and M is a right R module then we denote by Spec.R, the set of prime ideals of R and by min.Spec.R, the set of minimal prime ideals of R. Similarly Ass.M denotes the set of prime ideals of R associated to M on the right. Moreover r-ann.T denotes the right annihilator of a subset T of M and l-ann.T denotes the left annihilator of a left subset T of W in case W is a left R module. Also $|M|_r$ denotes the right Krull dimension of M if it exists and $|W|_l$ denotes the left Krull dimension of left R module W if it exists. For two subsets A and B of a given set A≤B means B contains A and A<B denotes A≤B but A≠B. Also A⊄B denotes the subset B does not contain subset A. For an ideal A of R c(A) denotes the set of elements of R that are regular modulo A

## §(1)  Main Theorem

We now prove results which culminate in the proof of theorem(4) which is our main theorem

Proposition (1): Let R be a right Noetherian, right Krull-homogenous ring with $|R|_r = \alpha$.

Let, $X = \{ P_i \in \text{spec.}R \;/\; |R/P_i|_r = |R|_r \}$.

Then there exists an ideal I of R such that:

(a) I is a nilpotent ideal of R and R/I is a right Krull-homogenous ring with $|R/I|_r = \alpha$ and

(b) If $P = \cap_{P_i \in X} P_i$; then I<P and P/I is a left localizable semi-prime ideal of the ring R/I

Proof: Fix an increasing critical composition of right ideals of R,

namely, $0 \neq A_1 < A_2 < A_3 < \ldots < A_{i-1} < A_i < \ldots < A_n = R$,

Where each factor module $A_i / A_{i-1}$ is a right $\alpha$- critical module.

Let $I_i = r\text{-ann.}(A_i / A_{i-1})$. Then, by [5], lemma(1.9), $R/I_i$ is a right $\alpha$ Krull-homogenous ring that is right $Q_i/I_i$-primary ring with

$Q_i/I_i = \text{Ass.}(A_i/A_{i-1})$. By [1], corollary (2.9), $R/I_i$ satisfies the right large condition. Hence there exists an ideal $J_i$ of R with $I_i < J_i$ and $J_i Q_i < I_i$. In fact, we can assume $J_i / I_i = l\text{-ann.}(Q_i/I_i)$. Thus $|R/j_i|_r < |R|_r = \alpha$. Note that if $I_i = Q_i$, then assume $J_i = R$. Now, put $J = \cap J_i$, $I = \cap I_i$ and $Q = \cap Q_i$.

Proof of (a): We now prove (a). First observe that since I ≤ r-ann.( $A_i/ A_{i-1}$); for $1 \leq i \leq n$, hence $I^n = 0$ and so I is nilpotent ideal of R. Next we show R/I is a right

α- K-homogenous ring . To seee this observe that if K/I is a right ideal of R/I with $|K/I|_r < |R/I|_r$, then for each j , $1 \leq j \leq n$ , we have that $(K+I_i)/I_j$ is a right ideal of $R/I_j$ such that $(K+I_i)/I_j \approx K/(I_j \cap K)$ which is a factor module of $K/I \cap K = K/I$ and hence

$| K/(I_j \cap K) |_r < |R/I|_r = |R/ I_j |_r$ . So $|(K+I_i)/I_j |_r < |R/ I_j |_r$ . But $R/ I_j$ is right α Krull-homogenous ring. Thus we must have $K \leq I_j$ , for each J , $1 \leq j \leq n$. Hence $K \leq \cap I_j = I$. Therefore R/I is right α K-homogenous ring with I a nilpotent ideal of R

Proof of (b): Next to prove (b) we first prove the following claim.

Claim: For each $P_i \varepsilon X$, there exists some $Q_i$ such that

$Q_i /I_i =$ Ass.( $A_i/ A_{i-1}$) and $P_i = Q_i$

Proof of the claim:- Observe first that JP ≤ I and its not difficult to see that $|R/J|_r < |R|_r$. Since by (a) above I is a nilpotent ideal with $I^n = 0$ so for each $P_i \varepsilon X$ we have that $0 \leq P_i$ implies that I ≤ P.

Also JQ ≤ I and $I \leq P_i$ implies $JQ \leq P_i$ . Obviously $J \not\leq P_i$ implies that Q ≤ $P_i$ Hence some $Q_j \leq P_i$ . Since each $P_i \varepsilon X$ is a minimal prime ideal of R , Hence $Q_j = P_i$. This proves the claim.

Hence Q= $\cap Q_i$ = P. Now JQ≤I implies since |R/J|<|R| , so P/I is a left localizable semi-prime ideal of R/I. This proves (b).

Lemma (2):- Let R be a right noetherian ring with $|R|_r = α$.

If an ideal T of R with $|R/T|_r = |R|_r = α$ is right weakly ideal invariant(w.i.i for short) , then for a finitely generated module M , if M has an essential α Krull-homogenous sub module H such that HT=0 and $|M/H|_r < |H|_r$, then MT=0.

Proof:- Given $|R|_r = \alpha$. Suppose first that M is a cyclic right R module, say, M=R/I, I a right ideal of R. Also we have that

$|M|_r = |H|_r = |R|_r = \alpha$. Now if H=L/I is an essential α Krull-homogenous right sub module of R/I, L a right ideal of R such that $|R/L|_r < |R|_r$; then if LT≤I, we will show that T≤I. To see this assume that T⊄I so that $|(T+I)/I|_r = |T/I \cap T|_r \leq |T/LT|_r < \alpha$, because T is right w.i.i. Clearly since L/I is essential in R/I, so L/I ∩ (T+I)/I ≠ 0. So there exists a non zero sub module A/I of L/I such that $|A/I|_r < |R|_r = \alpha$, a contradiction to L/I be α K-homogenous. Thus (T+I)/I =0. So T≤I. This proves the result for cyclic modules. Similarly we can prove this result for any finitely generated module. Hence MT=0.

Theorem (3):- Let R be a noetherian left and right Krull-homogenous ring with $|R|_r = \alpha$, and $|R|_l = \beta$. Let, X={ $P_i$ ε spec.R/ $|R/P_i|_r = |R|_r$ } and let $P = \cap_{P_i \varepsilon X} P_i$. Then $P^n = 0$, for some integer n≥1

Proof:- as in proposition (1) choose a nilpotent ideal I and an ideal J such that $|R/J|_r < \alpha$ and R/I is a right α-Krull homogenous noetherian ring. Moreover JP≤I. Suppose m is the nilpotent index of I. so $I^m = 0$ but $I^{m-1} \neq 0$. Clearly $I.I^{m-1} = 0$ implies

$JPI^{m-1} = 0$. Since $|R/J|_r < \alpha$, so by Goldie's theorem[3], there exists c ε c(P) ∩ J and $cPI^{m-1} = 0$. Thus if $PI^{m-1} \neq 0$ we can choose a non zero left ideal L of $PI^{m-1}$ such that QL = 0, where Q is a prime ideal of R such that $|R/Q|_l = \beta$. Clearly Q is a minimal prime ideal of R and I≤Q. Thus cL=0 implies $|L|_l < |R|_l = \beta$, a contradiction to R being left β-K-homogenous. Hence $PI^{m-1} = 0$. Again note that since I≤P, so $I^{m-1} < PI^{m-2}$ because $I^{m-1} = I.I^{m-2} \leq PI^{m-2}$. Now $PI^{m-1} = 0$ as above. Consider,

$A = PI^{m-2}/I^{m-1}$. Then since JP≤ I, and J≠0 we have that JA=0.

Now there exists c ε J ∩ c(P) as chosen above. Also observe that $I^{m-1}$ left large in $PI^{m-2}$, for otherwise there exists a non zero ideal $K \leq PI^{m-2}$ such that $K \cap I^{m-1} = 0$. Observe that $(K+I^{m-1})/I^{m-1} \approx K$ as left R-modules. Now JA=0 implies JK = 0. So cK=0, c ε c(P) as seen above. But R is left Krull homogenous so again using the same argument as above this also yields a contradiction. Hence $I^{m-1}$ is left large in $PI^{m-2}$. Also $I^{m-1}$ is left β-Krull-homogenous. Hence by lemma(2) above, we get that since $PI^{m-1} = 0$, so $P^2I^{m-2} = 0$. Continuing in this way we get a finite number of steps that $P^m = 0$

Theorem(4):- Let R be a noetherian ring which is left and right K-homogenous. Let $X = \{P_i \in \text{spec}.R / |R/P_i|_r = |R|_r\}$

and $P = \cap_{P_i \in X} P_i$.

Let $V = \{Q_j \in \text{spec}.R \mid |R/Q_i|_l < |R|_r\}$ and $Q = \cap_{Q_j \in V} Q_i$

Then the following hold true;

(1) $P^n = 0$, for some integer $n \geq 1$

(2) $|R/P_i|_r = |R|_r$ and $|R/P_i|_l = |R|_l$

(3) X=V

Proof:- (1) is true because of the previous theorem.

(2) Let $P_i$ be any prime ideal of X. By the left version of the previous theorem there exists some $m \geq 1$ such that $Q^m = 0$. So $Q^m \leq P_i$. Thus $Q \leq P_i$. Hence there exists some $Q_j$ in v such that $Q_j \leq P_i$. Since $P_i$ is a minimal prime ideal of R therefore $Q_j = P_i$. Thus $|R/P_i|_l = |R|_l$. This proves (2).

(3) By the argument in (2) above its clear that if $P_i \in X$ then $P_i \in V$. Thus $X \leq V$. Similarly we can show $V \leq X$. Therefore X = V and this proves (3).

References:

(1) A.K. Boyle, E.H. Feller; Semi critical modules and K-primitive rings, "L.M.N"; No. 700, Springer Verlag

(2) K.A. Brown, T.H. Lenagan and J.T. Stafford, "Weak Ideal Invariance and localisation"; J.Lond.math.society(2), 21 , 1980 , 53-61

(3) A.W. Goldie, "The Structure of Semi Prime Rings Under Ascending Chain Condition", Proc.Lond.math.Soc.8(1958)

(4) K.R. Goodaarl and R.B. Warfield , "An Introduction to Non commutative Noetherain Rings" , L.M.S., student texts ,16.

(5) R. Gordon , "some Aspects of Non Commutative Noetherian Rings" , L.N.M.; vol.545; Spg. Vlg. , 1975, 105-127